\documentclass[12pt]{amsart}
\usepackage{epsfig}
\usepackage{amsmath}
\usepackage{amsfonts}
\usepackage{amssymb}
\usepackage{mathabx}
\usepackage{epigraph, fancyhdr}

\def\C{\mathcal C}
\def\D{\mathcal D}
\def\F{\mathcal F}
\def\K{\mathcal K}

\def\L{\mathcal L} 
\def\H{\mathbf H}
\def\O{\mathcal O}

\def\1{\mathbf 1}
\def\M{\mathcal M}

\def\E{{\mathbf E}}

\def\ZZ{\mathbb Z}
\def\CC{\mathbb C}

\def\Res{\operatorname{Res}}

\def\hat{\widehat}
\def\tilde{\widetilde}
\def\p{\partial}
\def\a{\alpha}

\def\b{\beta}

\def\f{{\mathbf f}}
\def\g{{\mathbf g}}
\def\t{{\mathbf t}}

\def\ll{{\mathbf l}}
\def\vv{{\mathbf v}}

\def\gs{\sigma}

\def\k{\kappa}

\def\h{\hbar}

\def\lan{\langle}
\def\ran{\rangle}
\def\str{\operatorname{str}}

\def\ev{\operatorname{ev}}
\def\ft{\operatorname{ft}}

\def\td{\operatorname{td}}
\def\ch{\operatorname{ch}}

\def\eu{\operatorname{eu}}

\def\Eu{\operatorname{Eu}}

\def\tr{\operatorname{tr}}
\def\str{\operatorname{str}}

\def\square{\Box}
\def\und{\underline}
\renewcommand{\Delta}{\triangle}
\makeatletter
\DeclareFontFamily{OMX}{MnSymbolE}{}
\DeclareSymbolFont{MnLargeSymbols}{OMX}{MnSymbolE}{m}{n}
\SetSymbolFont{MnLargeSymbols}{bold}{OMX}{MnSymbolE}{b}{n}
\DeclareFontShape{OMX}{MnSymbolE}{m}{n}{
    <-6>  MnSymbolE5
   <6-7>  MnSymbolE6
   <7-8>  MnSymbolE7
   <8-9>  MnSymbolE8
   <9-10> MnSymbolE9
  <10-12> MnSymbolE10
  <12->   MnSymbolE12
}{}
\DeclareFontShape{OMX}{MnSymbolE}{b}{n}{
    <-6>  MnSymbolE-Bold5
   <6-7>  MnSymbolE-Bold6
   <7-8>  MnSymbolE-Bold7
   <8-9>  MnSymbolE-Bold8
   <9-10> MnSymbolE-Bold9
  <10-12> MnSymbolE-Bold10
  <12->   MnSymbolE-Bold12
}{}

\let\llangle\@undefined
\let\rrangle\@undefined
\DeclareMathDelimiter{\llan}{\mathopen}%
                     {MnLargeSymbols}{'164}{MnLargeSymbols}{'164}
\DeclareMathDelimiter{\rran}{\mathclose}%
                     {MnLargeSymbols}{'171}{MnLargeSymbols}{'171}
\makeatother

\title[Quantum Adams-Riemann-Roch]
      {Permutation-equivariant \\ quantum K-theory XI. \\
      Quantum Adams-Riemann-Roch}

\author[A. Givental]{Alexander GIVENTAL}
\thanks{This material is based upon work supported by the National 
Science Foundation under Grant  DMS-1611839, by the IBS Center for Geometry 
and Physics, POSTECH, Korea, and by IHES, France}

\begin{document}

\begin{abstract}

  We introduce {\em twisted} permutation-equivariant GW-invariants, and compute them in terms of untwisted ones. The computation is based on Grothendieck-like RR formula corresponding to Adams' operations from K-theory to itself, and the result can be understood as a ``quantum'' version of such Adams-RR. As in the case of cohomological quantum RR theorem \cite{CoGi}, the result is applied to express the invariants of bundle and super-bundle spaces in terms of those of the base. The bonus feature of permutation-equivariant K-theory is that the twisting classes can be understood as the simpler kappa-classes of Kabanov--Kimura \cite{KK}.     
  
\end{abstract}

\maketitle

\section*{Introduction}

In winter 1993-94, thinking of the Candelas {\em et al} mirror formula \cite{COGP} for quintic 3-folds, I arrived at the construction \cite{GiH} of the relevant hypergeometric functions from toric compactification of spaces of rational curves in the ambient projective $4$-space. About the same time I received an intriguing email from M. Kontsevich about his plan to obtain the numbers of degree-$d$ holomorphic spheres on the quintics by using Grothendieck-Riemann-Roch (GRR). It transpired soon that the starting point of his plan was the same as mine with toric compactifications, but applied to moduli spaces of stable maps instead. It has become standard: If $[\M]$ is the virtual fundamental class of such moduli space of degree-$d$ rational curves $\varphi: \Sigma \to\CC P^4$, then the virtual fundamental class of the moduli space for such curves in the quintic $X_5\subset \CC P^4$ is obtained from $[\M]$ by taking the cap-product with the Euler class of the vector bundle over $\M$ with the fibers $H^0(\Sigma; \varphi^*\O (5))$. The actual plan consisted in expressing the Euler class in terms of the components of the Chern character of this bundle, which can be computed using GRR and thus tracked back to some other GW-invariants of $\CC P^4$.

Kontsevich's plan was fully realized in 2001, in the form of the ``Quantum Lefschetz'' theorem based on ``Quantum Riemann-Roch'' \cite{CoGi}. Some fifteen years ago, discussing the subject with D. van Straten, I mentioned another idea: to use K-theory instead of cohomology, i.e. to count the number of degree-$d$ holomorphic spheres by the dimension of the virtual structure ring of the moduli space. His comment sounded intriguing too: ``Then you'll have to use Adams-Riemann-Roch.'' What he meant was that one would face the problem of expressing the K-theoretic analogue of the Euler class
of that same bundle in terms of Adams' operations, which are easier to compute following the GRR scheme.  

Namely, to a morphism between two abstract cohomology theories, there corresponds a Grothendieck-Riemann-Roch formula. The Adams-Riemann-Roch corresponds this way to Adams' operations $\Psi^k$, $k=\pm 1,\pm 2,\dots$, from K-theory to itself. Thus, a map $\pi: X\to Y$ between two compact complex manifolds induces the push-forward: $\pi_*: K^0(X)\to K^0(Y)$, which
intertwines with Adams' operations this way:
\[  \Psi^k(\pi_*(a)) = \pi_*\left( \Psi^k(a) \otimes \frac{\Eu (T_{X/Y})}{\Psi^k(\Eu (T_{X/Y}))}\right) .\]
Here $\Eu$ denotes the K-theoretic Euler class, defined on line bundles
by $\Eu(L)=1-L^{-1}$, so that the ratio, replacing the Todd class of the relative tangent bundle in the classical GRR, is characterized by
\[ \frac{\Eu(L)}{\Psi^k(\Eu(L))}=\frac{1-L^{-1}}{1-L^{-k}} .\]  

The present paper provides what we believe is the right answer to van Straten's question about ``Quantum Adams-Riemann-Roch''. The previous approaches included the Quantum Hirzebruch-RR \cite{Co} by T. Coates for {\em fake} K-theoretic (and, even more generally, cobordism-valued) GW-invariants (see also \cite{CGL, GiF}), and the adelic characterization \cite{ToN} by V. Tonita of twisted K-theoretic GW-invariants in genus $0$. It turns out that the subject fits naturally into the setting of permutation-equivariant quantum K-theory. Among the applications we include the abstract analogues of ``Quantum Serre'' and ``Quantum Lefschetz'' theorems. A detailed discussion of more concrete applications, involving K-theoretic mirror formulas and $q$-hypergeometric functions, is postponed to a sequel paper. 

\section{Twistings}

Let $X$ be a compact K\"ahler manifold. As usual, we denote by $X_{g,n,d}$ the moduli space of degree-$d$ stable maps to $X$ of nodal compact connected complex curves of arithmetical genus $g$ carrying $n$ marked points.
Consider the diagram
\[ \begin{array}{ccc} X_{g,n+1,d} & \stackrel{\ev}{\longrightarrow} & X \\
  \ft \downarrow  &  & \\
  X_{g,n,d} & &  \end{array}, \]
where $\ev$ and $\ft$ denote evaluation at and forgetting of the last, $n+1$-st marked point. The diagram can be viewed as the degree-$d$ {\em universal stable map} to $X$ of genus-$g$ complex connected nodal curves with $n$ marked points. Namely, when $\varphi: (\Sigma, \sigma_1,\dots, \sigma_n) \to X$ is such a stable map, the fiber of $\ft$ over its equivalence class $[\varphi]\in X_{g,n,d}$ is canonically identified with $\Sigma /\operatorname{Aut}(\varphi)$, and the restriction of $\ev$ to this fiber coincides with the map defined by $\varphi$. 

Given a vector bundle $E \in K^0(X)$ over $X$, we put $E_{g,n,d}:=\ft_*\ev^*E$. Here $\ft_*$ denotes the full K-theoretic push-forward.\footnote{Think of the virtual bundle $H^0(\Sigma, \varphi^*E)\ominus H^1(\Sigma, \varphi^*E)$.}
It is not hard to show (see e.g. \cite{CoGi}) that the orbisheaf $E_{g,n,d}$ has a locally free resolution and thus represents an element in the Grothendieck group $K^0(X_{g,n,d})$ of orbibunldes over $X_{g,n,d}$.

In cohomological GW-theory, one defines twisted GW-invariants of $X$ by systematically replacing virtual fundamental cycles $[X_{g,n,d}]$ by their
cap-products with a general invertible multiplicative characteristic class of $E_{g,n,d}$. Mimicking (and generalizing) this construction, we define {\em twisted} K-theoretic GW-invariants of $X$. Namely, for each non-zero integer $k$, fix an element $E^{(k)} \in K:=K^0(X)\otimes \Lambda$, and replace the virtual structure sheaf $O_{X_{g,n,d}}$
introduced in \cite{YPLee} with the tensor product
\[ \O^{\E^{(\bullet)}}_{X_{g,n,d}}:= \O_{X_{g,n,d}} \otimes e^{\textstyle \sum_{k\neq 0} \Psi^k(E^{(k)}_{g,n,d})/k},\]
where we use the superscript $\E^{(\bullet)}$ to refer cumulatively to the twisting datum $\{ E^{(k)} \}$ (which can be considered as a $(\ZZ-0)$-graded bundle). Here $\Lambda$ is a ground ring, which is a $\lambda$-algebra, i.e. is equipped with the action of Adams' operations $\Psi^k$. For the sake of applications, and in contrast with previous papers of this series, we assume
that they are defined on $\Lambda$ for negative values of $k$ as well, and defined the action of such operations on Novikov's variables by $\Psi^{-1}(Q^d)=Q^d$. The factors $1/k$ are placed in the exponent for convenience, having in mind the applications where half of all $E^{(k)}$ (say, with $k<0$, or with $k>0$) are equal to each other, while the other half are zeroes, resulting in the twisting factors of the form $\Eu (E_{g,n,d})$ or $1/\Eu (E^*_{g,n,d})$. In such applications, it is necessary to assume that the K-theory in question is equivariant with respect to the fiberwise scalar action of the circle $\CC^{\times}$ on the bundles $E^{(k)}$. The same trick is useful in order to battle potential divergence issues. For example, to guarantee convergence of the expression for $\O^{tw}_{X_{g,n,d}}$, one can extend the ground ring $\Lambda$ into $\Lambda [\lambda_{+}^{\pm 1}, \lambda_{-}^{\pm 1}]$ by multiplicative coordinates $\lambda_{\pm}$ on two circles, and assume that each $E^{(k)}$ is a multiple of $\lambda_{+}$ for $k>0$ and  of $\lambda_{-}$ for $k<0$. Since $\Psi^k(\lambda_{\pm})=\lambda_{\pm}^k$, the expression for $\O^{tw}_{X_{g,n,d}}$ comes out as a power series in $\lambda_{+}$ and $1/\lambda_{-}$ with well-defined coefficients.

In parallel with the construction of Part IX, we introduce the total descendant potential $\D_X^{tw}$ of {\em twisted} permutation-equivariant quantum K-theory of $X$. 
It is a function of a sequence $\t = (\t_1,\t_2,\dots)$ of $K$-valued Laurent polynomials in one variable, $q$ (which is a place-holder for universal cotangent line bundles $L_i$), and depends on the choice of $E^{(k)}\in K$ considered as parameters. In this notation, 
\[ \D_X^{tw}: =e^{\textstyle \sum_{r>0} \h^{(g-1)r} \Psi^r(\F_g^{\E^{(r\bullet)}}(\t_r,\t_{2r},\dots))/r}, \]
where
\[ \F_g^{\E^{(r\bullet)}}(\t_r,\t_{2r},\dots) :=
\sum_{d, \ll} \frac{Q^d}{\prod_r l_r!}\ \lan \t_r(L),\dots ; \t_{2r}(L),\dots; \dots\ran_{g,\ll,d}^{\E^{(r\bullet)}} .\]
We refer the reader to Part IX for the definition of the correlators, but  remind that $\ll = (l_1,l_2,\dots)$ is a partition of the number $n=\sum_r r l_r$ of the marked points (into $l_r$ cycles of length $r$, permuted by the elements of $S_n$). 

Note the superscripts $\E^{(r\bullet)}$, indicating that the twisted structure sheaves participating in the definition of the correlators and generating functions are based on the twisting sequence $E^{(rk)}$, $k=\pm 1,\pm 2,\dots$, depending on $r$.

\medskip

\section{A digression}

We need to make a digression in order to justify the formula for $\D_X^{tw}$.

Let us recall that the actual definition of the total descendant potential in permutation-equivariant quantum K-theory has the structure
\[ \sum_{\eu} \h^{-\eu} \sum_{m=0}^{\infty} \frac{1}{m!} \sum_{g \in S_n} \str_g H^*\left(\M_{\eu, m} ; V_{\eu,m} \right) ,\]
where $\M_{\eu,m}$ are the moduli spaces of $m$-component stable maps of Euler characteristic $\eu$, and $V_{\eu,m}$ are appropriate $S_m$-equivariant orbisheaves. This can be rewritten according to the cycle structure of permutations:
\[ \exp \left\{ \sum_{\eu} \h^{-\eu} \sum_{r=1}^{\infty} \frac{1}{r} \str_h H^*(\M_{\eu,r} ; V_{\eu,r}) \right\} ,\]
where $h$ is induced by the cyclic permutation of $r$ connected components of the curves. More explicitly, the terms contributing to each $\str_h$ have the following form:
\[ \str_h H^*\left( \M^r ; \O_{\M^r}\otimes e^{\sum_{k\neq 0} \Psi^k(\sum_{i=1}^r pr_i^*E^{(k)}_{g,n,d})/k} \prod_{i=1}^r pr_i^*(T) \right),\]
where $\M$ is the nickname for $X_{g,n,d}$, $pr_i: \M^r\to \M$ is the projection to the $i$th factor, $T$ is the tentative notation for the tensor product over $\M$ of all inputs from the marked points, and the action of the generator $h \in \ZZ_r$ on the sheaf cohomology is induced by the cyclic permutation of the factors. We claim that $\str_h$ coincides with
\[ \Psi^r \left( H^* \left(\M; \O_{\M}\otimes e^{\sum_{k\neq 0} \Psi^k(E_{g,n,d}^{(rk)})/k} T \right)\right) ,\]
implying the above formula for $\D_X^{tw}$. The essence of our claim is captured by the following abstract lemma. 

\medskip

{\tt Lemma.} {\em Let $\M$ be a compact complex manifold (or, more generally, virtual orbifold), $E, T \in K^0(\M)\otimes \Lambda$, $\pi:\M \to pt$, and $F$ a polynomial (or convergent series) in one variable with $\Psi$-invariant coefficients. Then   
  \begin{align*} \str_h (\pi^r)_* \left( F\left(\frac{\Psi^k(pr_1^*E+\cdots +pr_r^*E)}{k}\right)\otimes \prod_{i=1}^r pr_i^*T\right)  \\
    = \left\{ \begin{array}{ll}\Psi^r\left( \pi_* \left(F\left(\frac{\textstyle \Psi^{k/r}(E)}{\textstyle k/r}\right)\otimes T\right)\right) & \text{when $r \mid k$}  \\
      \Psi^r (\pi_*(F(0) \otimes T)) & \text{when $r \nmid k$} \end{array} \right. \end{align*}
}

{\tt Proof.} Due to Lefschetz-Kawasaki's localization formula (see Part IX),
$\str_h (\pi^r)_*(W)$ is determined by $\tr_h (W|_{\M})$, where $\M$ is the fixed point locus of $h$ in $\M^r$. The restriction of
$\sum_i pr_i^* E$ is $E\otimes \CC[\ZZ_r]$, where the 2nd factor carries the regular representation of $\ZZ_r$. We have: $\tr_h \Psi^k(\CC[\ZZ_r]) = 0$ when $r\nmid k$, and $=r$ when $r\mid k$. In the latter case, we end up with
\[ \tr_h \left( F\left(\frac{\Psi^k(E\otimes \CC[\ZZ_r])}{k}\right)\otimes T^{\otimes r} \right) =  \Psi^r \left( F\left(\frac{\Psi^{k/r}(E)}{k/r}\right)\otimes T\right) .\]
This coincides with $\str_h (W|_{\M})$ when
\[ W=\prod_{r=1}^r pr_i^*\left( F\left(\frac{\Psi^{k/r}(E)}{k/r}\right)\otimes T\right) .\]
Now $\str_h (\pi^r)_*(W)$ can be computed explicitly and yields
\[ \Psi^r \left(\pi_*\left(F\left(\frac{\Psi^{k/r}(E)}{k/r}\right)\otimes T\right) \right).\]
The simpler case $r\nmid k$ is processed similarly.

\medskip

\section{Formulations}

Let us return now to our main problem of expressing the family $\D_X^{tw}$ in terms of $\D_X$.

In order to formulate the answer, we interpret the total descendant potential as a quantum state $\lan \D_X^{tw}\ran$ in the Fock space associated to a certain symplectic loop space $(\K^{\infty}, \Omega^{\infty})$. By definition,
$\K^{\infty}$ consists of sequences $\f = (\f_1,\f_2,\dots)$ or $K$-valued rational functions of $q$ with possible poles at $q=0,\infty$, or at roots
of unity. The $\Lambda$-valued symplectic form is defined by
\[ \Omega^{\infty}(\f,\g):=\sum_{r>0} \frac{\Psi^r}{r}(\Omega^{(r)} (\f_r,\g_r)),\] where 
\[ \Omega^{(r)} (f, g) :=-[ \Res_{q=0}+\Res_{q=\infty} ] \ (f(q^{-1}), g(q))_{\E^{(r\bullet)}}\ \frac{dq}{q},\]
and $(\cdot,\cdot)_{\E^{(r\bullet)}}$ is the {\em twisted} K-theoretic Poincar\'e pairing in $K^0(X)\otimes \Lambda$:
\[  (a,b)_{\E^{(r\bullet)}}:= \chi \left(X; a \otimes b \otimes e^{\textstyle \sum_{k\neq 0} \Psi^k(E^{(rk)})/k} \right) .\]
Thus, the symplectic structure $\Omega^{\infty}$ depends in fact on parameters $E^{(k)}$.

The loop space $\K^{\infty}$ is equipped with polarization $\K^{\infty}=\K^{\infty}_{+}\oplus \K^{\infty}_{-}$ which is Lagrangian with respect
to $\Omega^{\infty}$ for any values of the parameters. Namely, $\K_{+}$ consists
of sequences $\f=(\f_1,\f_2,\dots)$ of $K$-valued Laurent polynomials in $q$,
while $\K_{-}$ consists of sequences of rational functions such that $\f_r(\infty)=0$ and $\f_r(0)\neq \infty$ for each $r=1,2,\dots$. The function
$\D_X^{tw}$ is naturally defined on $\K^{\infty}_{+}$. More precisely, we assume
that the argument $\t = (\t_1,\t_2,\dots)$ consists of Laurent polynomials
$\t_r$ whose coefficients are {\em small} in $K=K^0(X)\otimes \Lambda$ in the sense that they belong to some ideal in $K$, filtration by whose powers is increased by Adams' operations $\Psi^r$ with $r>0$. We define $\lan \D_X^{tw}\ran$ as an element (more precisely, a one-dimensional subspace) in the Fock space by shifting the argument $\t$ by the
{\em dilaton vector} $\vv = ((1-q)\1, (1-q)\1, \dots )$, where $\1$ is the unit in $K^0(X)$, and then lifting the dilaton-shifted function from $\K^{\infty}_{+}$
to $\K^{\infty}$ by the projection $\K^{\infty} \to \K_{+}^{\infty}$ along $\K^{\infty}_{-}$. Thus, $\lan \D_X^{tw}\ran$ is a family of functions (in a
family Fock spaces) defined in a small neighborhood of the dilaton vector.   

\medskip

{\tt Theorem 1} (Quantum Adams-Riemann-Roch).
{\em  $\lan \D_X^{tw} \ran = \widehat{\square} \lan \D_X \ran$, where $\widehat{\square}$ is quantization of the operator $\square: \K^{\infty} \to \K^{\infty}$  acting on $\f =(\f_1,\f_2,\dots)$ as component-wise multiplication
  $\f_r \mapsto \square_r(q) \f_r$, where
\[ \square_r(q) = e^{\textstyle \sum_{k\neq 0} \Psi^k(E^{(rk)})/k(1-q^{k})} .\]}    
\indent We remind that our rules of quantization for quadratic hamiltonians have the form
\[ \widehat{q_\a q_\b} =\h^{-r} q_{\a}q_{\b}, \ \widehat{q_{\a}p_{\b}} =q_{\a}\p_{q_{\b}}, \ \widehat{p_\a p_\b}= \h^r \p_{q_{\a}}\p_{q_{\b}},\]
where $ p_{\a}, q_{\b} $ form a set of Darboux coordinates on the component
$\K^{(r)}$ of the symplectic loop space $\K^{\infty}=\prod_{r>0} \K^{(r)}$
equipped with the symplectic form $\Psi^r\Omega^{(r)}/r$. This results in a well-defined {\em projective} representation of the Lie algebra of quadratic hamiltonians. Symplectic transformations are quantized by $\hat{M}=e^{\widehat{ \log M}}$. Note that due to the identity $(1-q^k)^{-1}+(1-q^{-k})^{-1}=1$, the operator of multiplication by $\square_r$ satisfies
\[ \square_r (q^{-1})\ \square_r (q) = e^{\textstyle \sum_{k\neq 0} \Psi^k(E^{(rk)})/k} .\]
Therefore it is not symplectic on $\K^{(r)}$, but rather represents a family of operators (parameterized by $\E^{(r\bullet)}$) which transform the symplectic form $\Omega^{(r)}$ based on the twisted pairing $(\cdot,\cdot)_{\E^{(r\bullet)}}$ into the undeformed symplectic form $\Omega^{(r)}|_{\E=0}$ based on the Poincar\'e pairing $(a,b)=\chi(X; a\otimes b)$. Respectively, the quantization $\widehat{\square}$ trivializes the family of Fock spaces by mapping their elements (which are interpreted as functions on the underlying symplectic spaces) in the opposite direction: from the untwisted Fock space to the twisted ones. Thus, Theorem 1 says that under this trivialization, the family of quantum states $\lan \D_X^{tw}\ran$ remains constant.

\medskip

Our entire setup can be generalized by allowing on the role of the twisting elements $E^{(k)}$ Laurent polynomials of the form $E=\sum_m E_m q^m$ with coefficients $E_m\in K$, and by defining the bundles $E^{(k)}_{g,n,d}$ as
\[ E_{g,n,d}:=\ft_*\left( \sum_{m\in \ZZ} \ev^*(E_m) \ L^m \right).\]
Here $L$ is the {\em relative} universal cotangent line on the universal curve diagram, i.e. $L=L_{n+1}$ on $X_{g,n+1.d}$. The following theorem shows that
when all $E^{(k)}$ are divisible by $1-q$, the effect of such twisting on the generating function $\D_X$ consists in a shift of the origin. 

\medskip

{\tt Theorem 2.} {\em The total descendant potential $\D_X^{tw}$ of permutation-equivariant quantum K-theory, twisted by the sequence $E^{(k)}$, $k\neq 0$, of $K$-valued Laurent polynomials in $q$, is related to  
the quantum state $\widehat{\square}\ \lan \D_X \ran$, where $\square$ is defined by componentwise multiplication operators (from Theorem 1)
  \[ \f_r\mapsto \square_r (q) \f_r  = e^{\textstyle \sum_{k\neq 0} \Psi^k(E^{(k)}(1))/k(1-q^k)} \f_r \]
corresponding to the sequence of twisting elements $E^{(k)}|_{q=1}$, by the change of the dilaton vector $\vv = (\vv_1,\vv_2,\dots)$ from $\vv_r=(1-q)\, \1$ for all $r=1,2,\dots$ to
\[ \vv_r = (1-q)\, e^{\textstyle \sum_{k\neq 0} \Psi^k (E^{(kr)}(q)-E^{(kr)}(1))/k(1-q^k)}.\]}
\indent We remind that Adams' operations act on $q$ by $\Psi^k(q)=q^k$, and hence
$\Psi^k(E(q)-E(1))/(1-q^k)=\Psi^k[(E(q)-E(1))/(1-q)]$ is a $K$-valued Laurent polynomial in $q$.

Note that the components $\vv_r$ of the dilaton vector here depend on $r$. This
should be understood this way: as the parameter $E(q)-E(1)$ varies, the quantum state $\lan D_X^{tw}\ran$ stays unchanged, but the construction of this quantum state from the family of generating functions $\D_X^{tw}$ for GW-invariants involves the varying dilaton shift.  

\section{Remarks on kappa-classes}

In the context of intersection theory in Deligne-Mumford spaces
$\overline{\M}_{g,n}$, one defines $\kappa_{m-1}:=\ft_*(\psi_{n+1}^m)$, where
$\ft: \overline{\M}_{g,n+1}\to \overline{\M}_{g,n}$ is the universal family of curves, and $\psi_{n+1}=c_1(L_{n+1})$ is the 1st Chern class of the relative cotangent line bundle.

In \cite{KK}, Kabanov and Kimura generalize this notion in the context of GW-theory by calling kappa-classes the push-forwards by $\ft=\ft_{n+1}: X_{g,n+1,d}\to X_{g,n,d}$ of any class of the form $\sum_{m\geq 0}\ev_{n+1}^*(\phi_m) \psi_{n+1}^m$, and prove a simple formula describing the effect of introducing  into GW-intersection theory all kappa-classes with $\phi_0=0$. Namely, a general dilaton shift ensues.

In fact Theorem 2 can be viewed as the complete analogue of Kabanov -- Kimura's formula in the case of permutation-equivariant quantum K-theory. Indeed, given a K-valued Laurent polynomial $E=\sum_m E_m q^m$, the K-theoretic push-forward
\[ E_{g,n,d}:=\ft_* \left(\sum_{m\in \ZZ} \ev^*(E_m) L^m\right)\]
defines the K-theoretic version of a kappa-class in the sense of Kabanov-Kimura. In the permutation-equivariant context, one is interested, in particular, in the traces of cyclic permutations acting on the tensor products of $k$ copies of $E_{g,n,d}$, which yield $\Psi^k(E_{g,n,d})$. In complete analogy with Kabanov-Kimura's result, the effect of kappa-classes with $E(1)=0$ on the total descendant potential is described by Theorem 2 as the dilaton shift, while the case when $E(1)\neq 0$ is also covered. 

\section{Quantum Serre}

Let $\D_X^{\E}$ and $\D_X^{\E^*}$ denote the total descendant potentials of permutation - equivariant quantum K-theory on $X$, twisted in the first case by a sequence $\E:=\{ E^{(k)} \}$ of bundles $E^{(k)}\in K$, and in the second by the sequence $\E^*=\{ \Psi^{-1}(E^{(-k)})\} $ of dual bundles indexed in the reverse order. Thus, for a fixed $k$, in the former case we use $e^{\Psi^k(E_{g,n,d})/k}$ as the twisting factor, and $e^{-\Psi^k(E_{g,n,d})^*/k}$ in the latter. The name of the result comes from the fact that $E_{g,n,d}$ and $-(E_{g,n,d})^*$ are related by Serre's duality. 

\medskip

{\tt Corollary 1.} 
{\em The symplectic identification  $(\K^{\infty},\Omega^{\infty}_{\E}) \to (\K^{\infty},\Omega^{\infty}_{\E^*})$ given by the operators of multiplication $\f_r \mapsto e^{\sum_{k\neq 0} \Psi^k(E^{(rk)})/k}\f_r$ transforms the quantum state $\lan \D_{X}^{\E^*}\ran$ into $\lan \D_X^{\E} \ran$.}

\medskip

Indeed, according to Theorem 1,
\begin{align*} \lan \D_X^{\E}\ran &= \widehat{\square}_{\E} \lan \D_X\ran, \ \text{where}\ \log (\square_{\E})_r =  \sum_{k\neq 0} \frac{\Psi^{k}(E^{(rk)})}{k(1-q^k)},\\
 \lan \D_X^{\E^*} \ran &= \widehat{\square}_{\E^*} \lan \D_X \ran, \ \text{where}\ 
 \log (\square_{\E^*})_r = \sum_{k\neq 0} \frac{\Psi^k(\Psi^{-1}E^{(-rk)})}{k(1-q^k)}.\end{align*}
Thus $\lan \D_X^{\E} \ran = \widehat{\square}_{\E} \widehat{\square}_{\E^*}^{-1} \lan \D_X^{\E^*} \ran$, where
$(\square_{\E})_r(\square_{\E*})_r^{-1} = e^{\sum_{k\neq 0} \Psi^k(E^{(rk)})/k}.$

\medskip

{\tt Remarks.} (1) The corollary does not mean that $\D_X^{\E}$ and $\D_X^{\E^*}$ practically coincide as generating functions for GW-invariants. Indeed, the corresponding inputs $\t_r^{\E}$ and $\t_r^{\E^*}$ of these generating functions can be found from
\[ e^{\sum_{k\neq 0} \Psi^k(E^{(kr)})/k} \left((1-q)\1+\t_r^{\E}\right) = (1-q)\1+\t_r^{\E^*}.\]
\indent (2) Quantum Serre duality can be easily extended to the twistings given by sequences $E^{(k)}$ of $K$-valued Laurent polynomials in $q$. In this case, the identification of the quantum states involves multiplications by $e^{\sum_{k\neq 0} \Psi^k(E^{(rk)}(1))/k}$, but in the comparison of the generating functions, the additional dilaton shifts occur. 

\section{Eulerian twistings}

As it was explained in Introduction, the interest in studying twisted GW-invariants is fueled by applications to submanifolds $Y\subset X$ given by holomorphic sections of a bundle $E$ over $X$. This requires the twisting by the Euler classes $\Eu (E_{g,n,d})$, but to make them invertible, one first considers a $\CC^{\times}$-equivariant version of the theory (making $\lambda\in \CC^{\times}$ act on $E$ by fiberwise scalar multiplication) and then passes to the limit $\lambda \to 1$ when possible. The approach actually works only in genus $0$, but the Euler-twisted GW-theory can be considered in any genus, and we call the resulting GW-theory the quantum K-theory of the {\em super-manifold} $\Pi E$ (which is motivated by $\operatorname{sdim} \Pi E = \dim Y$).

For a line bundle $L$, we have
\[ \Eu (L)=1-L^{-1} = e^{\textstyle \sum_{k<0} \Psi^k(L)/k},\]
and hence by the splitting principle
\[ \Eu (E_{g,n,d}) = e^{\textstyle \sum_{k<0} \Psi^k(\ft_*(E))/k}.\]
Thus, the total descendant potential $\D_{\Pi E}$ of the suprmanifold $\Pi E$ is obtained from the twisting sequence $E^{(k)}=E$ for $k<0$ and $E^{(k)}=0$ for $k>0$.

\medskip

For another application, consider K-theoretic GW-invariants of a non-compact space defined as the total space of a vector bundle, $F$, over $X$. In comparison with the previous case, it is convenient to think of $F$
as dual to $E$, equipped with the dual action of the circle, i.e. $\lambda\in \CC^{\times}$ acts on $F$ fiberwise as multiplication by $\lambda^{-1}$. The invariants of $F$ can be defined by fixed point localization: $X_{g,n,d}$ are the fixed point loci of the $\CC^{\times}$-action
on moduli space of curves in the bundle space, while the twisting bundle, coming from the {\em denominator} $\Eu (F_{g,n,d})$ in Lefschetz' fixed point localization formula, is
\[ \Eu^{-1}(F_{g,n,d})=e^{\textstyle -\sum_{k<0} \Psi^k(\ft_*(F))/k}.\]
Thus, the total potential $\D_F$ of the bundle space is obtained from the twisting sequence $E^{(k)}=-F=-E^*$ for $k<0$ and $0$ for $k>0$.

\medskip

{\tt Corollary 2.}{\em
  \begin{align*} \lan \D_{\Pi E}\ran &= \left[ e^{\textstyle \sum_{k<0} \Psi^k(E)/k(1-q^k)} \right]^{\widehat{\ \ }} \lan \D_X \ran , \\
   \lan \D_{E^*}\ran &= \left[ e^{\textstyle -\sum_{k<0} \Psi^k(E^*)/k(1-q^k)} \right]^{\widehat{\ \ }}\lan \D_X \ran .\end{align*} }

Indeed, the operators $\square_r$ from Theorem 1 in either Eulerian case do not depend on $r$.

\medskip

{\tt Remark.} Note that the multiplication operators in this corollary
are asymptotical expansions of the following infinite products:
\begin{align*} \prod_{l=0}^{\infty} \Eu (Eq^{l}) \sim\ & e^{\sum_{k<0} \Psi^k(E)/k(1-q^k)} \sim \prod_{l=1}^{\infty} \Eu^{-1}(E q^{-l}) \\
  \prod_{l=0}^{\infty} \Eu^{-1} (E^* q^{l}) \sim\ & e^{\sum_{k<0} -\Psi^k(E^*)/k(1-q^k)} \sim \prod_{l=1}^{\infty} \Eu (E^* q^{-l}) .\end{align*}
Namely, for a line bundle $E$, we have the formal expansion:
\[ \prod_{l=0}^{\infty}\Eu(Lq^l)=\prod_{l=0}^{\infty}(1-L^{-1}q^{-l}) \sim e^{\sum_{l=0}^{\infty} \sum_{k<0} L^kq^{kl}/k} =e^{\sum_{k<0} L^k/k(1-q^k)}. \]

\section{Genus 0}

In Part X, from the all-genus adelic description of $\D_X$, we derived the adelic characterization of the range $\L_X\subset \K^{\infty}$ of the ``big J-function'' constructed from the dilaton-shifted genus-0 descendant potential $\F_0$ of permutation-equivariant quantum K-theory on $X$.

The theory carries over almost {\em verbatim} to our setting of twisted GW-invariants. Let $\E^{(\bullet)}=\{ E^{(k)} \}$ be the twisting sequence consisting of $K$-valued Laurent polynomials in $q$. Denote by $\L_X^{tw}$ the subvariety in $\K^{\infty}$ (depending on $\E^{(\bullet)}$) which consists of sequences $\f=(\f_1,\f_2,\dots)$ of $K$-valued rational functions in $q$ of the form
\begin{align*} & \f_r = (1-q)e^{\sum_{k\neq 0}\Psi^{k}(E^{(rk)}(q)-E^{(rk)}(1))/k(1-q^k)}+\t_r(q)+ \\ & \sum_{\ll,d,\a}\frac{Q^d e^{-\sum_{\k\neq 0} \Psi^k(E^{(rk)}(1))/k}\phi_\a}{\prod_{s>0} l_s!} \lan \frac{\phi^{\a}}{1-qL}, \t_r(L),\dots ; \t_{2r}(L),\dots ; \dots \ran_{0,\ll+\1_{1}\,d}^{\E^{(r\bullet)}},\end{align*}
where $\t = (\t_1,\t_2,\dots) \in \K^{\infty}_{+}$. Here $\{ \phi_\a \}$ and $\{ \phi^\a \}$ are bases in $K$ Poincar\'e-dual with respect to the undeformed pairing. The superscript $\E^{(r\bullet)}$ indicates that the correlators are defined using the virtual structure sheaves twisted by $e^{\sum_{k\neq 0} \Psi^k(E^{(rk)}_{0,n+1,d})/k}$, where $n=1+\sum_{r>0} rl_r$, and the partition $\ll+\1_1$ consists of $l_1+1$ cycles of length $1$ and $l_r$ cycles of lengths $r>1$. 

\medskip

{\tt Corollary 3} (quantum Adams-RR in genus 0): {\em
$\square \, \L_X^{tw} = \L_X$,
where $\square$ acts on $\f=(\f_1,\f_2,\dots)\in \K^{\infty}$ by $(\square\, \f)_r=e^{ \sum_{k\neq 0} \Psi^k (E^{(rk)}(1))/k(1-q^k)}\ \f_r$.}

\medskip  

As we will see, this result needs a comment rather than a proof. Let
\[ \F_0^{\E^{(\bullet)}}(\vv+\t):=\sum_{0,\ll,d} \frac{Q^d}{\prod_{s>0}l_s!} \lan \t_1(L),\dots; \t_2(L)\dots; \dots \ran_{0,\ll,d}^{\E^{(\bullet)}}\]
be the dilaton-shifted version of the twisted genus-0 descendant potential, where $\vv_r$ are as in Theorem 2 (and the same as in the top line of the above formula for $\f_r$). One can describe $\L_X^{tw}$ in terms of the differentials of $\F_0^{\E^{(r\bullet)}} (\t_r,\t_{2r},\dots)$ considered as families of functions of $\t_r$ depending on $\t_{2r},\t_{3r},\dots$ as parameters:
$\f_r(\t_r,\t_{2r},\dots) = \t_r+d_{\t_r}\F_0^{\E^{(r\bullet)}}$.
The way how $\F_0^{\E^{(\bullet)}}$ is obtained from its untwisted value $\F_0:=\F_0^{\E^{(\bullet)}}|_{\E=0}$ is described by Theorem 2. Namely, $\lan \D_X^{tw}\ran = e^{\textstyle \F_0^{\E^{(\bullet)}}/\h +\dots}$, where the ellipsis denotes terms of other orders in $\h$. Thus, it suffices to figure out how the quantized operator $\widehat{\square}$ acts on $\lan \D_X \ran$, and extract the terms of order $h^{-1}$ from $\log (\widehat{\square} \lan \D_X \ran)$.

In fact the multiplication operator $\square_r$ has one important property which simplifies the description of its quantization: it preserves the negative space $\K^{(r)}_{-}$ of the Lagrangian polarization on $(\K^{(r)}, \Omega^{(r)})$. This implies that, in Darboux coordinate notation, the quadratic hamiltonian of $\log \square_r$ consists of monomials $p_\a q_\b$ and $q_a q_\b$, but does not contain $p_\a p_\b$. In other words, $\widehat{\square}_r$
act by a linear change of variables and multiplication. An explicit description is given in \cite{GiQ} (see Proposition 5.3). According to it
\[ \widehat{\square}_r\ \D (\t_r) = e^{\textstyle W_r(\t_r)/\h^r} \D ([\square_r^{-1} \t_r]_{+}),\]
where $[\dots ]_{+}$ is the projection to $\K_{+}^{(r)}$ along $\K_{-}^{(r)}$,
and $W_r$ is a certain quadratic form in $\t_r$, whose explicit description in terms of the operator $\square_r$ we presently omit. Extracting the terms of
order $\h^{-1}$, we find
\[ \F_0^{\E^{(\bullet)}}(\t_1,\dots, \t_r, \dots) = W_1(\t_1) + \F_0
([\square_1^{-1} \t_1]_{+}, \dots, [\square_k^{-1} \t_k]_{+},\dots) .\]
Consequently
\[  W_r(\t_r) + \F_0([\square_r^{-1} \t_r]_{+}, \dots, [\square_{rk}^{-1} \t_{rk}]_{+},\dots) =\F_0^{\E^{(r\bullet)}}(\t_r, \dots, \t_{rk},\dots) .\]
This shows that, indeed, $\f=(\f_1,\f_2,\dots)$ lies in $\L_X^{tw}$ whenever 
$\square \f = (\square_1\f_1,\square_2\f_2,\dots)$ lies in $\in \L_X$. 

\section{Quantum Lefschetz}

Continuing our study of genus-$0$ theory, consider the Eulerian cases, and denote by $\L_X$, $\L_{\Pi E}$, and $\L_{E^*}$ the varieties in the appropriate loop spaces representing the quantum K-theory of $X$, of the super-manifold $\Pi E$, and of the total space of the bundle $E^*$. We will have to assume that $E$ is the direct sum of line bundles $E_i$, $i=1,\dots, s$, or at least that Chern roots of $E$ are integer.

\medskip

{\tt Theorem 3} (Quantum Lefschetz). {\em 
  Let $\f = (\f_1,\f_2,\dots)$ be a point in $\L_X$, represented by power series in Novikov's variables $\f_r = \sum_d \f_{r,d} Q^d$ with coefficients $\f_{r,d}$ which are $Q$-independent vector-valued rational functions of $q$. Then $\f^{\Pi E}\in \L_{\Pi E}$ and $\f^{E^*}\in \L_{E^*}$, where
\begin{align*} \f^{\Pi E}_r &:= \sum_d f_{r,d}Q^d \prod_{i=1}^s \frac{\prod_{l\leq (c_1(E_i),d)} (1-E_i^{-1}q^l)}{\prod_{l\leq 0} (1-E_i^{-1}q^l)} , \\
  \f^{E^*}_r &:= \sum_d f_{r,d}Q^d \prod_{i=1}^s \frac{\prod_{l<(c_1(E_i),d)} (1-E_iq^{-l})}{\prod_{l<0} (1-E_iq^{-l})} .\end{align*}
}

{\tt Proof.} Recall that in Part X, we established $\D_q$-invariance properties of $\L_X$. Let $P_\a q^{Q_\a\p_{Q_\a}}$ be the translation operators which, together with the operators of multiplication by $Q_\a$, generate $\D_q$, and
let $E_i = \prod_\a P_\a^{-m_{i\a}}$ be the representation of $E_i$ in terms of the line bundles $P_\a$. By definition, we have $(c_1(P_\a), d) = -d_\a$, and hence
\[ D_i(d):=(c_1(E_i), d) = \sum_{\a} m_{i\a} d_{\a} .\]
Consider the finite difference operators obtained by the asymptotical expansions of the following infinite products:
\begin{align*} \Gamma_{\Pi E} &\sim \prod_{i=1}^s \prod_{l\leq 0} \left(1-q^l \prod_\a (P_\a q^{Q_\a\p_{Q_\a}})^{m_{i\a}}\right) , \\
  \Gamma_{E^*} &\sim \prod_{i=1}^s \prod_{l\leq 0} \left(1-q^l \prod_\a (P_\a q^{Q_\a\p_{Q_\a}})^{-m_{i\a}} \right)^{-1} . \end{align*}
Then we have:
\begin{align*} \Gamma_{\Pi E} Q^d  &\sim Q^d \prod_{i=1}^s \prod_{l=-\infty}^ {D_i(d)} (1-E_i^{-1}q^l) \sim Q^d \square_{\Pi E} \prod_{i=1}^s \frac{\prod_{l=-\infty}^{D_i(d)} (1-E_i^{-1}q^l)}{ \prod_{l=-\infty}^0 (1-E_i^{-1}q^l)} , \\
 \Gamma_{E^*} Q^d  &\sim Q^d \prod_{i=1}^s \prod_{l=-\infty}^{-D_i(d)} (1-E_i^{-1}q^l)^{-1} \sim Q^d\square_{E^*}\prod_{i=1}^s\frac{\prod_{l=-\infty}^0(1-E_iq^l)}{\prod_{l=-\infty}^{-D_i(d)}(1-E_iq^l)} .\end{align*} 
By the $\D_q$-invariance property of $\L_X$, if $\f\in \L_X$, then $\Gamma_{\Pi E}\f \in \L_X$ and $\Gamma_{E^*}\f\in \L_X$. By the genus-0 Adams-RR, this implies that $\square_{\Pi E}^{-1} \Gamma_{\Pi E} \f$ lies in $\L_{\Pi E}$, and $\square_{E^*}^{-1} \Gamma_{E^*} \f$ lies in $\L_{E^*}$. It remains to examine each $Q^d$-term to conclude that $\square_{\Pi E}^{-1} \Gamma_{\Pi E} \f = \f^{\Pi E}$, and
likewise $\square_{E^*}^{-1} \Gamma_{E^*} \f = \f^{E^*}$ since
\[ \frac{\prod_{l=-\infty}^0(1-E_iq^l)}{\prod_{l=-\infty}^{-D_i(d)}(1-E_iq^l)} = \frac{\prod_{l<(c_1(E_i),d)} (1-E_iq^{-l})}{\prod_{l<0} (1-E_iq^{-l})} .\]

{\tt Remark.} The $\D_q$-invariance properties of $\L_X$ extend to $\L_{\Pi E}$ and $\L_{E^*}$, e.g. because according to the genus-0 Adams-RR, $\L_{\Pi E}$ and $\L_{E^*}$ are related to $\L_X$ by the multiplication operators \[ \square_{\Pi E}=e^{\textstyle \sum_{k<0}\Psi^k(E)/k(1-q^k)}\ \text{and}\  \square_{E^*}=e^{\textstyle -\sum_{k<0}\Psi^k(E^*)/k(1-q^k)}\]
respectively, which commute with finite difference operators from $\D_q$. Alternatively, one can check that the whole theory of adelic characterization of $\L_X$ from Part X carries over to the twisted case.

\section{An example} 

Let $P$ denote the Hopf bundle over $X=\CC P^{n-1}=\operatorname{proj}(\CC ^n)$. It generates $K^0(\CC P^{n-1})$ and satisfies the relation $(1-P)^n=0$. As can be shown by elementary arguments (see \cite{GiL}), the ``small J-function'' in quantum K-theory of $\CC P^{n-1}$ has the form
\[ J:=(1-q)\sum_{d\geq 0}\frac{Q^d}{\prod_{l=1}^d(1-Pq^l)^n}.\]
The projection $[J]_{+}=1-q$ of $J$ to $\K_{+}$  shows that $J$ corresponds to the input $\t=0$. This implies that in permutation-equivariant quantum K-theory of $\CC P^{n-1}$, the sequence $\f = (\f_1,\f_2,\dots)$ all of whose components $\f_r=J$ represents a point in $\L_{\CC P^{n-1}}$ corresponding to the input $\t=(\t_1,\t_2,\dots)$ with all $\t_r=0$.

The cotangent bundle\footnote{I am thankful to P. Koroteev for discussions of quantum K-theory of quiver varieties \cite{KPSZ}, whetting my interest in the cotangent bundles.} of $\CC P^{n-1}$ can be described as $Hom (\CC^n/P, P)$, i.e. as $T^*=nP-1$ in $K^0_T(\CC P^{n-1})$. By $-1$, the subtraction of the trivial 1-dimensional bundle $Hom (P,P)$ is meant. Our general theory does not prohibit twistings by virtual sums of line bundles, but to remove the degeneration of the twisted Poincar\'e pairing, we need to equip the whole bundle with a non-trivial action of $\CC^{\times}$, whose elements we will denote by $\lambda$. Now applying the quantum Lefschetz theorem, we find a point  $\f^{T^*\CC P^{n-1}}$ in $\L_{T^*\CC P^{n-1}}$, all of whose components are equal to
\[ I:=(1-q)\sum_{d=1}^{\infty}Q^d \frac{\prod_{l=0}^{d-1}(1-\lambda P^{-1}q^{-l})^n}{\prod_{l=1}^d(1-Pq^l)^n}.\]
The subtraction of the trivial bundle in $nP-1$ does not have any
effect on the formula for $I$. However, it contributes the factor $(1-\lambda)$
to the equivariant Euler class
\[ \Eu (T^*\CC P^{n-1}) = \frac{(1-P^{-1}\lambda)^n}{(1-\lambda)} ,\]
and hence in the twisted Poincar\'e pairing:
\[ (a, b)_{T^*\CC P^{n-1}} = \chi \left(\CC P^{n-1}; \frac{(a\otimes b) (1-\lambda)}{(1-P^{-1}\lambda)^n}\right) .\]
Note that the factor $1/(1-P^{-1}\lambda)^n$, which has no non-equivariant limit
since at $\lambda=1$, the denominator vanishes in $K^0(\CC P^{n-1})$, is present in the numerator of each term in $I$ with $d>0$. Therefore, when components of $I$ are extracted by pairing with a basis (say, $\phi_\a =(1-P)^\a$, $\a=0,\dots,n-1$), all terms with $d>0$ will have non-equivariant limits.
This is expected, since moduli spaces of positive degree stable maps to $T^*\CC P^{n-1}$ are compact. Indeed, for a non-constant map $\varphi: \Sigma \to \CC P^{n-1}$ to the zero section of the cotangen bundle, the exact sequence $0\to T^*\CC P^{n-1}\to O(-1)\otimes \CC^n \to \O \to 0$ on $\CC P^{n-1}$ induces
the exact cohomological sequence
\begin{align*}  0&\to H^0(\Sigma;\varphi^*T^*\CC P^{n-1}) \to H^0(\Sigma; \varphi^*\O(-1)\otimes \CC^n) \to H^0(\Sigma;\O) \to \\
  &\to  H^1(\Sigma,\varphi^* T^*\CC P^{n-1}) \to H^1(\Sigma, \varphi^*\O(-1)\otimes \CC^n) \to H^1(\Sigma;\O) \to 0,\end{align*}
where $H^0(\Sigma; \varphi^*\O(-1))=0$, and hence $H^0(\Sigma; \varphi^*T^*\CC P^{n-1})=0$. Too bad,  though, that due to the factor $1-\lambda$ the non-equivariant limit of $(I-(1-q))/\Eu (T^*\CC P^{n-1})$ vanishes identically.

The situation changes\footnote{Which is easy to realize from the above long exact sequence (where $H^1(\Sigma; \O)=0$ too when $g=0$).}, if one consders the twisting bundle $E^*=nP = \O(-1)\otimes \CC^n$. The series $\f_r^{E^*}$ turn out to be the same as $I$, but the formula for the Poincar\'e pairing does not contain the factor $1-\lambda$ this time, and as a result the part with $d>0$ has a meaningful (and rather simple) non-equivariant limit:
\[ \lim_{\lambda\to 1} \frac{I-(1-q)}{\Eu(E^*)} = (1-q)\sum_{d>0} \frac{(-Pq^{d/2})^{-n(d-1)}\, Q^d}{(1-Pq^d)^n} .\] 
\indent Note that the projection $[I]_{+}$ of $I$ to $\K_{+}$ is rather complicated. Thus,  $\f = (I,I,\dots)$  represents a point in $\L_{E^*}$ and in $\L^{T^*\CC P^{n-1}}$ corresponding to a non-trivial input $\t$.
Nonetheless, as it follows from the ``explicit reconstruction'' results of Part X, starting from the series $I$ (respectively $J$), one can obtain an explicit parameterization of the entire varieties $\L_{E^*}$ and $\L_{T^*\CC P^{n-1}}$ (and respectively $\L_{\CC P^{n-1}}$).
   
\section{Hodge bundles}   

When $E$ is the $1$-dimensional trivial bundle, $E_{g,n,d}=1-\H_g^*$, where $\H_g$ is the {\em Hodge bundle} of dimension $g$. Its fibers are formed by holomorphic differentials on the Riemann surfaces, and do not depend on the numbers of marked points, or the maps to the target space. The quantum Adams-RR theorem says that the twisting of the virtual structure sheaves
\[ \O_{X_{g,n,d}} \mapsto  \O_{X_{g,n,d}}\otimes
e^{\textstyle  \sum_{k\neq 0} \Psi^k(\tau_k (1-\H_g^*))/k}, \ \ \tau_k\in \Lambda_{+},\]
results in the transformation of the total descendant potential
\[ \lan \D_X \ran \mapsto \widehat{\square}\, \lan \D_X\ran,\ \text{where} \ \square_r=e^{\textstyle \sum_{k\neq 0}\Psi^k(\tau_{rk})/k(1-q^k)}, r=1,2,\dots\]
In the case $g=0$, the Hodge bundle is zero, and the twisting of the virtual structure sheaves by the scalar $e^{\sum_{k\neq 0} \Psi^k(\tau_k)/k}$ results in the multiplication of the whole genus-0 descendant potential $\F_0(\t_1,\t_2,\dots)$ by this scalar. Furthermore, the genus-$0$ part of $\log D_X$, which has the form $\sum_{r>0} \h^{-r} \Psi^r (\F_0 (\t_r,\t_{2r}, \dots))/r$, turns into
\[ \sum_{r>0} \h^{-r} \Psi^r(\Delta_r \F_0(\t_r,\t_{2r},\dots))/r, \ \text{where} \ \Delta_r=e^{\textstyle \sum_{k\neq 0} \Psi^k(\tau_{rk})/k}.\] 
On the other hand, on the loop space $\K^{\infty}$, the symplectic form $\Omega^{\infty}(\f,\g)$ $= \sum_{r>0} \Psi^r \Omega (\f_r,\g_r)/r$ turns into its  twisted form
\[ \Omega^{\infty}_{tw}(\f,\g)=\sum_{r>0} \Psi^r ( \Delta_r \Omega (\f_r,\g_r) )/r.\] The factors
$\Delta_r$ occur here because of the change of the standard Poincar\'e pairing 
into $(a,b)_{\E^{(r\bullet)}} = \chi (X; a\otimes b\otimes \Delta_r)$. Note that the components of $\Omega^{\infty}$ and of the genus-0 part of $\log \D_X$ change the same way under the twisting. Recalling the construction of $\L_X^{tw} \subset \K^{\infty}$, we conclude that the Hodge twisting does not affect it at all. Now the genus-0 quantum Adams RR implies that {\em $\L_X$ is invariant under the operators $\square: \f = (\f_1,\f_2,\dots)\mapsto (\square_1 \f_1, \square_2 \f_2,\dots)$.} In Part X, we derived the same result under the name ``string flows'' from adelic characterization of $\L_X$ as a special case of $\D_q$-invariance properties.  

\section{Proof of Theorem 1}

The main difficulty in direct application of Adams--Grothendieck-RR technique to our problem is caused not by the permutations of the marked points, but by their interaction with the orbifold structure of the moduli spaces (due to the intrinsic automorphisms of stable maps). Yet, Theorems 1 and 2 look as if such stacky issues could be ignored. It is possible therefore that our results can be derived by direct K-theoretic arguments similar to those used by Kabanov-Kimura \cite{KK} or by Tonita \cite{ToT} in the case of cohomological kappa-twistings. Lacking such a direct argument, we resort to another approach: the ``quantum Hirzebruch-RR'' theorem of Part IX, which reduces the problem to fake quantum K-theory, where the stacky issues can be ignored indeed.  

\medskip

In Part IX, we gave an adelic representation of $\D_X$ in terms of
$ \D_{X/\ZZ_M}^{tw}$, the generating functions for twisted in a certain way fake K-theoretic GW-invariants of the orbifold target spaces $X/\ZZ_M=X\times BZ_M$, $M=1,2,\dots$. The representation was based on the virtual Kawasaki-RR formula \cite{ToK} applied to computing the traces of the transformations in sheaf cohomology on $X_{g,n,d}$ induced by renumbering of the marked points. Let us recall that Kawasaki's formula computes such traces in terms of certain fake holomorphic Euler characteristics on strata of the inertia orbifold, which in our case were described in terms of moduli spaces of stable maps to $X/\ZZ_M$. More precisely, the sum of the contributions of
all the strata has the structure of Wick's summation over decorated graphs, where the vertex contributions (i.e. contributions of one-vertex graph without edges) have the form of certain twisted fake K-theoretic GW-invariants of $X/\ZZ_M$. 

The whole machinery and the resulting description fully applies\footnote{In fact, there is some change in the formalism; namely the Poincar\'e pairing $(\cdot,\cdot)^{tw}_r$ depends on $r$, while in the untwisted case it is the same for all $r$, but as we will see shortly, this does not produce any destructive effect.} to the case of $\D_X^{tw}$. That is, vertex contributions
now come in the form of $\D_{X/\ZZ_M}^{tw}$ where however the change in the structure sheaves $\O_{X_{g,n,d}} \mapsto \O_{X_{g,nd}}^{\E^{(\bullet)}}$ results in the additional twistings of the virtual fundamental classes in moduli spaces of stable maps to $X/\ZZ_M$. Since we don't have any more room for decorating $\D_{X/\ZZ_M}^{tw}$, we will use this notation for such additionally twisted fake K-theoretic GW-invariants of $X/\ZZ_M$, and will use the more cumbersome symbol $\D_{X/\ZZ_M}^{tw}|_{\E=0}$ to denote the generating function used in Part IX in the description of $\D_X = \D_X^{tw}|_{\E=0}$.

In order to specify the additional twisting, let us denote by $\M$ a moduli space of stable maps to $X/\ZZ_M$, considered as a Kawasaki stratum in some $X_{g,n,d}$, and recall from Part IX the general shape of its contribution
to $\str_h H^*(X_{g,n,d}; \O_{X_{g,n,d}}\otimes V)$, where $h$ is some transformation on $X_{g,n,d}$ defined by a renumbering of marked points, and $V$ is an $h$-equivariant bundle on $X_{g,n,d}$.

A point in $\M$ is represented by a map $\varphi: \hat{\Sigma} \to X$ of a compact holomorphic curve to $X$, equipped with a (possibly ramified, and not necessarily connected) principal $\ZZ_M$-cover $p: \Sigma\to \hat{\Sigma}$. The composition $p \circ \varphi$ represents a stable map to $X$ fixed under the
symmetry $h$ defined by the generator in $\ZZ_M$. The contribution of the Kawasaki stratum has the form
\[  \chi^{fake}\left( \M; \frac{\tr_h (V|_{\M})}{\str_h \wedge^{\bullet} N^*_{\M}}\right),\]
where $N^*_{\M}$ is the virtual conormal bundle of $\M$ considered as a stratum in a moduli space of covering stable maps $p\circ \varphi$ to $X$, and $\wedge^{\bullet}$ denotes the exterior algebra. By definition,
\[ \chi^{fake} (\M ; W) := \int_{[\M]} \td (T_{\M}) \ \ch (W),\]
where $[\M]$ is the virtual fundamental cycle, and $T_{\M}$ is the virtual tangent bundle. We also remind that
$\tr_h W$ is defined on an $h$-invariant vector bundle $W$ as $\sum_{\lambda} \lambda W_{\lambda}$, where $W_{\lambda}$ is the eigen-bundle corresponding to the eigenvalue $\lambda$.  

In the untwisted theory, the bundle $V$ comes from inputs into the correlators
of quantum K-theory on $X$ which the formula purports to compute. The denominators provide the twisting of the fake quantum K-theory of $X/\ZZ_M$, which results in replacing $\D_{X/\ZZ_M}^{fake}$ with what we currently denote $\D_{X/\ZZ_M}^{tw}|_{E=0}$. In our current situation, $V$ also contains the twisting factor $e^{\sum_{k\neq 0} \Psi^k(E^{(k)}_{g,n,d})/k}|_{\M}$. Thus, we have to replace $\chi^{fake}(\M; W)$ with
\[ \chi^{tw}(\M; W) := \chi^{fake} \left(\M: W\otimes e^{\textstyle \sum_{k\neq 0} \tr_h \Psi^k(E^{(k)}_{g,n,d}|_{\M})/k}\right) .\]
We claim that this modification can be described as the appropriate twisting of the fake quantum K-theory on $X/\ZZ_M$ in the same sense in which this term is used in \cite{ToT} and earlier in this paper.

The trick is the same as the one used in Part IX (and elsewhere). For each $M$th root of unity $\lambda = e^{2\pi i a/M}$, $a=1,\dots, M$, denote by $\CC_{\lambda}$ the $1$-dimensional complex orbibundle on $X/\ZZ_M$ on which the generator $h$ of $\ZZ_M$ acts with the eigenvalue $\lambda$. Let $\ft: \C \to \M$ denotes the universal curve over $\M$ (considered as a moduli space of stable maps to $X/\ZZ_M$), and $\ev: \C \to X/\ZZ_M$ denotes the universal stable map in this sense. Then $\ft_* (\ev^*(E\otimes \CC_{\lambda^{-1}}))$ is identified with the eigen-subbundle in $E^{(k)}_{g,n,d}|_{\M}$ on which $h$ acts with the eigenvalue $\lambda$. After the application of $\Psi^k$, the eigenvalue will become $\lambda^k$. Thus, $\chi^{tw}(\M ; W) =$ 
\[ \chi^{fake} \left( \M ; W \otimes \prod_{a=1}^M
e^{\textstyle \sum_{k\neq 0} e^{2\pi i a k/M} \Psi^k[\ft_*\ev^*(E^{(k)}\otimes \CC_{e^{-2\pi i a/M}})]/k}\right) .\]

Now we are in the right position to invoke the result of V. Tonita \cite{ToT} (actually going back to the ``orbifold quantum Riemann-Roch theorem'' of H.-H. Tseng \cite{Ts}) about several superimposed twistings of GW-invariants of orbifold target spaces. It is formulated in cohomological terms, but since fake K-theoretic invariants are defined purely cohomologically, the rephrasing in K-theoretic notations is straightforward.   

In fact the computation is similar to what we've already done in Part IX. 
The family of quantum states $\lan \D_{X/\ZZ}^{tw}\ran$ in the appropriate symplectic loop space formalism is obtained from the state at the value $\E=0$ of the parameters by quantized multiplication operators acting separately in each {\em sector} labeled by elements $h^s\in \ZZ_M$, $s=1,\dots, M$. Combining the twisting for all values of $a=1,\dots, M$, but for a fixed value of $k$, we find the multiplication operator in the form
\[ S_M^{(k)}(E^{(k)}) \prod_{a=1}^{M} \prod_{l=1}^{\infty} S_a^{(k)}(q^{l-\{as/M\}} E^{(k)}), \text{where}\ S_a^{(k)} (\cdot) = e^{\textstyle e^{2\pi i ka/M} \Psi^k(\cdot )/k},\]
and $\{ x/M \}$ denotes the fractional part. Just as in Part IX, let $r=(s,M)$ be the greatest common divisor, $m=M/r$, $s'=s/r$, so that $(s',m)=1$, and write $a=a't'+mu$, where $t'$ is inverse to $s'$ modulo $m$, $0\leq a'<m$, and $u=1,\dots, r$. Then $\{ as/M \} = a'/m$ regardless of the value of $u$. Collecting the terms in the exponent with fixed values of $k, a', s'$, we take into account that $\sum_{u=1}^r e^{2\pi i ku/r}=0$ unless $k$ is divisible by $r$, in which case it is equal to $r$. For $k=rk'$, we
also write $e^{2\pi i k a't'/M} = \zeta^{a'k'}$, where $\zeta =e^{2\pi i t'/m}$ is a {\em primitive} $m$th root of unity. Thus, the exponent in our infinite product becomes the product of $\Psi^{rk'}(E^{(rk')})/k'$ with the geometric series
\[ 1+\sum_{a'=0}^{m-1} \sum_{l=1}^{\infty} \zeta^{k'a'} q^{rk'(l-a'/m)}) = \sum_{n=0}^{\infty} q^{rk'n/m}\zeta^{-k'n} 
  = \frac{1}{1-(q^{r/m}/\zeta)^{k'}} .\]
  Renaming $k'$ by $k$ again, and summing over all $k\neq 0$, we can summarize this computation by saying that in the sector $h^s$ determined by the {\em level} $r=(s,M)$ and the primitive $m$th root of unity $\zeta$ (such that $\zeta^{t'}=e^{2\pi i/m}$ where $t'$ is inverse to $s'=s/r$ modulo $m$) the multiplication operator has the form
 \[ \square_r^{(\zeta)}= \Psi^r \left( e^{\textstyle \sum_{k<0} \Psi^k(E^{(kr)})/k(1-(q^{1/m}/\zeta)^{k})} \right) .\]
Notice that in notations of Theorem 1, \  
$\square_r^{(\zeta)} (q) = \Psi^r \left( \, \square_r (q^{1/m}/\zeta) \right)$.

\medskip
 
We now combine these results for all $M=1,2,\dots$ and interpret them   in terms of the adelic symplectic loop space and its quantization.

As we've already said, $\D^{tw}_{X/\ZZ_M}$, after dilaton shift by $(1-q^M)\1$ in the unit sector, is interpreted as a quantum state $\lan \D_{X/\ZZ_M}^{tw}\ran$
in a suitable symplectic loop space (depending on parameter $\E$) equipped with a certain Lagrangian polarization. The family of quantum states is obtained form $\lan \D_{X/\ZZ}^{tw}|_{\E=0}\ran$ by quantization of the operator
acting as multiplication by $\square_r^{(\zeta)}$ sector-wise.

In more detail, the sector labeled by the level $r$ and a root of unity $\zeta$ of primitive order $m=m(\zeta)$ (so that $M=mr$) is represented in the symplectic loop space by a copy $\K_r^{(\zeta)}$ of $\K^{fake}=K((q-1))$.
The $\Lambda$-valued symplectic form pairs $\K_r^{(\zeta)}$ with $\K_r^{(\zeta^{-1})}$ according to the rule 
\[  \K_r^{(\zeta)}\times \K_r^{(\zeta^{-1})} \ni (f,g) \mapsto \Omega_r^{(\zeta)}:=
\frac{1}{rm(\zeta)} \Res_{q=1} (f(q^{-1}), g(q))^{(r)} \, \frac{dq}{q}.\]
Here $(a,b)^{(r)}:= r \Psi^r (\Psi^{1/r}(a),\Psi^{1/r}(b))^{tw}_r$, that is, more explicitly
\[ (a,b)^{(r)} = r \Psi^r\left( \chi \left(X; \Psi^{1/r}(ab)\, e^{\textstyle \sum_{k\neq 0}\Psi^k(E^{(rk)})/k}\right) \right).\]
Note that $\Omega_r^{(\zeta)}(f,g)=\Omega_r^{(\zeta)}\left(\square_r^{(\zeta)}f, \square_r^{(\zeta^{-1})}g\right)|_{\E=0}$, since
\[ \square_r^{(\zeta)}(q^{-1}) \square_r^{(\zeta^{-1})}(q) = \Psi^r \left(e^{\textstyle \sum_{k\neq 0} \Psi^k(E^{(rk)})/k} \right).\]
This means that the operators $\square_r^{(\zeta)}$ taken together are not symplectic but rather identify the symplectic form depending on $\E$ with its
value at $\E=0$. At the quantum level, these operators transform $\lan \D^{tw}_{X/\ZZ_M}|_{\E=0}\ran$ into $\lan \D^{tw}_{X/\ZZ_M}\ran$. 

We should add here that the dilaton-shifted functions $\D_{X/\ZZ_M}^{tw}$, which depend, beside Novikov's variables $Q$, on the Planck constant $\h$, and the sectorial inputs $t_r^{(\zeta)}\in K[[q-1]]$ with $m(\zeta)r=M$, are homogeneous in such a way that
$\lan \D_{X/\ZZ_M}^{tw} \ran (t, \h, Q) = \lan \D_{X/\ZZ_M}^{tw}\ran (t/\sqrt{\h}, 1, Q)$. Indeed, this is true for $\E=0$, and remains true after the application of the quantized operators $\square_r^{(\zeta)}$ due to the homogeneity of our quantization rules.

\medskip

Now, we combine together all the sectors $\K_r^{(\zeta)}$ for all $M=m(\zeta)r$
to form the {\em adelic} symplectic loop space $\und{\K}^{\infty}=\prod_{\zeta}\prod_r \K_r^{(\zeta)}$. It is equipped with the symplectic form $\und{\Omega}^{\infty} = \oplus_{r,\zeta} \Omega_r^{(\zeta)}$, and
with Lagrangian polarization $\und{\K}^{\infty}_{\pm}$ (which comes from the splitting of each $\K_r^{(\zeta)}$ into the direct sum of subspaces $(\K_r^{(\zeta)})_{\pm}$ described in Part IX), and the dilaton vector
$\und{\vv} = (\und{\vv}_1, \und{\vv}_2,\dots)$ with $\und{\vv}_r=(1-q^r)\1 \in \K_r^{(\1)}$. Using these structure, we declare the adelic tensor product
\[ \lan \und{\D}_X^{tw}\ran (\{ t_r^{(\zeta)} \}, \h, Q) := \bigotimes_{M=1}\left(\lan \D_{X/\ZZ_M}^{tw}\ran (\{ t_r^{(\zeta)}/\sqrt{\h}\}_{rm(\zeta)=M}, 1, Q^M)\right) ,\]
which is naturally a function on $\und{\K}^{\infty}_{+}$, a quantum state 
in the quantization of the adelic symplectic loop space $(\und{\K}^{\infty},\und{\Omega}^{\infty})$. Our previous discussion based on the application of the ``orbifold quantum Riemann-Roch theorem'' can be summarized as
\[  \lan \und{\D}_X^{tw} \ran = \widehat{\und{\square}} \, \lan \und{\D}_X^{tw}|_{\E=0} \ran ,\]
where the adelic operator $\und{\square}$ acts block-diagonally as sector-wise multiplication by $\square_r^{(\zeta)}$.   

\medskip

Now we recall from Part IX how Wick's summation over graphs leads to the adelic description of $\lan \D_X^{tw} \ran$. The {\em adelic map}: $\und{\ }: \K^{\infty} \to \und{\K}^{\infty}$ associates to a sequence $\f = (\f_1,\f_2,\dots)$ a collection
$\und{\f} = \{ f_r^{(\zeta)} \}$ of $K$-valued Laurent series in $q-1$, where $f_r^{(\zeta)}$ is the Laurent series expansion near $q=1$ of $\Psi^r(\f_r(q^{1/m}/\zeta)$. It is straightforward to check that the adelic map is symplectic, respects the dilaton vectors, and maps $\K^{\infty}_{+}$ to $\und{\K}^{\infty}_{+}$, but does not respect the negative spaces of the polarizations. Moreover, the adelic image of $\K^{\infty}_{-}$ forms together with $\und{\K}_{+}^{\infty}$ a new polarization in $\und{\K}^{\infty}$, which we call {\em uniform}, while the initial polarization $\und{\K}_{\pm}^{\infty}$ is called {\em standard}.\footnote{For the sake of completeness, we reproduce here its definition. Let $\phi_{\a}$ and $\phi^{\a}$ run Poincar\'e-dual bases in $K$. Then $\Psi^r \left(\phi_{\a}(q^{1/m)}-1)^n\right)$, $n=0,1,2,\dots$, run a basis in $(\K_r^{(\zeta)})_{+}:=K[[q-1]]$ (assuming that $m=m(\zeta)$), while the Darboux-dual basis in $(\K_r^{(\zeta^{-1})})_{-}$ is run by $\Psi^r\left( (\phi^\a q^{n/m}e^{\textstyle -\sum_{k\neq 0} \Psi^k(E_{rk})/k})/(1-q^{1/m})^{n+1} \right) $.}
It turns out that the effect of the edge propagators in Wick's summation over graphs on the vertex tensor product $\lan \und{\D}_X^{tw}\ran $ is equivalent to representing the same quantum state in the uniform polarization (instead of the standard one). The main result of Part IX, extended to our current (twisted) situation, describes the quantum state $\lan \D_X^{tw}\ran$ as {\em induced} by the adelic map from $\lan \und{\D}_X^{tw}\ran$, i.e. as simply the restriction of the function to the adelic image of $\K^{\infty}_{+}$ in $\und{\K}^{\infty}_{+}$ {\em after} the polarization's change.

\medskip

We are ready now to derive Theorem 1 from the adelic description of $\D_X^{tw}$
and our previous computation of $\und{\D}_X^{tw}$. Namely, our observation
$\square_r^{(\zeta)}(q)=\Psi^r(\square_r(q^{1/m}/\zeta))$ shows that the adelic map intertwines the operator $\square$ on $\K^{\infty}$ with $\und{\square}$ on
$\und{\K}^{\infty}$, i.e. for any $\f\in \K^{\infty}$, the adelic image of $\square\ \f$ coincides with $\und{\square} \ \und{\f}$. In order to conclude from this that
\[ \widehat{\square} \ \lan \D_X^{tw}|_{\E=0}\ran = \text{pull-back of}\  \widehat{\und{\square}}\ \lan \und{\D}_X^{tw}|_{\E=0}\ran \ \text{by adelic map} = \lan \D_X^{tw}\ran , \]
we need to make sure that, applying the quantized operator $\widehat{\und{\square}}$ to $\lan \und{\D}_X^{tw}|_{\E=0}\ran$, and then restricting the outcome to the adelic image of $\K^{\infty}_{+}\subset \und{\K}^{\infty}_{+}$, we get the same result as when restricting first and then applying $\widehat{\square}$ to $\D_X^{tw}|_{\E=0}$.

\medskip

This can be derived from the property of the expressions
\[ \square_r(q) = e^{\textstyle \sum_{k\neq 0} \Psi^k(E_{rk})/k(1-q^{k})}\]
which we have already exploited in the section of genus-0.   
The exponent is an $\operatorname{End}(K)$-valued rational function of $q$ with poles at the roots of unity, but no poles at $q=0$ or $\infty$. This shows that the operator $\square$ preserves $\K_{-}^{\infty}$. The adelic extension $\und{\square}$ of the operator $\square$ preserves the subspace $\K^{\infty}$ embedded by the adelic map into $\und{\K}^{\infty}$ and preserves the negative space of the uniform polarization --- because it coincides with the image of $\K^{\infty}_{-}$. This is indeed a specific feature of our infinite dimensional situation: the same subspace is Lagrangian in both the symplectic subspace, and the ambient symplectic space. In naive notation of Darboux coordinates $(p_\a, q_\a)$, a quadratic Hamiltonian whose flow preserves the Lagrangian subspace given by the equations $q_\a=0$ for all $\a$, must be free of terms $p_\a p_\b$ (i.e. contain only monomials of the form $q_\a p_\b$ or $q_\a q_\b$). Under quantization, $q_\a q_\b$ becomes the the operator of multiplication by $q_{\a}q_{\b}/\h$, which commutes with the restriction of functions on a subspace. The bilinear terms $\sum_{\a\b} q_\a A_{\a\b} p_\b$ under quantization yield $\sum_{\a\b} q_{\a}A_{\a\b} \p_{q_{\b}}$ which is the vector filed corresponding to the (half of the) Hamilton equations: $\dot{q}_\b = \sum_\a q_{\a}A_{\a\b}$. The property of this vector field to be tangent to a certain subspace ($\K^{\infty}_{+}\subset \und{\K}^{\infty}_{+}$ in our situation) is equivalent to the property of the derivation to commute
with the operation of restriction of functions to this subspace. This observation completes the proof of Theorem 1.

\medskip

{\tt Remark.} The last argument, though convincing, isn't really necessary.
In other words, if the exponent contained an operator of multiplication by a rational function of $q$ (and hence not preserved $\K_{-}^{\infty}$) its quantization would still commute with the operation of restriction of functions from $\und{\K}^{\infty}_{+}$ to $\K_{+}^{\infty}$.  
The reason is that the adelic image of $\K^{\infty}$ though does not coincide with the entire $\und{\K}^{\infty}$, is dense in some week topology, and anyway behaves as
a coisotropic subspace. Namely, exact sequence
\[ 0\to \K^{\infty} \to \und{\K}^{\infty} \to \und{\K}_{+}^{\infty}/\K_{+}^{\infty}\to 0 \]
shows that linear functions $L$ defining the subspace $\K^{\infty}$ inside $\und{\K}^{\infty}$ become operators of multiplication under quantization.  
Let $H$ be a quadratic hamiltonian whose flow preserves the subspace.
Then the Poisson bracket $\{ H, L\}$ is another function vanishing on the subspace. Under quantization, Poisson brackets between quadratic and linear hamiltonians are preserved: $[\widehat{H},\widehat{L}]=\widehat{ \{ H,L\} }$. Consequently, for an element $L \D$ in the Fock divisible by $L$, we find that  $\widehat{H}L\D=L\widehat{H}\D+\{H,L\}\D$ also vanishes on the subspace.

 \section{Proof of Theorem 2}

Following the same logic as in our proof of Theorem 1, we arrive at the twisted fake quantum K-theory of $X/\ZZ_M$ with the additional twisting of the ``fake virtual structure sheaf'' of a moduli space $\M$ by
 \[ \prod_{a=1}^M
 e^{\textstyle \sum_{k\neq 0} e^{2\pi i a k/M} \Psi^k[\ft_*\ev^*((E^{(k)}(L)-E^{(k)}(1))\otimes \CC_{e^{-2\pi i a/M}})]/k} ,\]
where $E^{(k)}$ is a sequence of $K$-valued Laurent polynomials in $q$. 
This time we apply another ``quantum Riemann--Roch theorem'' of V. Tonita \cite{ToT} (about twistings of type B). This result describes the effect of such twisting as an additional dilaton shift {\em in the unit sector}. Namely, the dilaton shift, defining the quantum state $\lan \D_{X/\ZZ_M}^{tw}\ran$ with zero twisting bundles (as in Part IX), as well as with the twisting bundles $E^{(k)}(1)$, has the form $(1-q^M)\1$. According to the result of \cite{ToT} it should be replaced with
\[ (1-q^M) \prod_{k\neq 0} \prod_{a=1}^M e^{\textstyle e^{2\pi i ak/M} \Psi^k[(E^{(k)}(q)-E^{(k)}(1))/k(1-q)]}.\]
We see that the exponent for a fixed $k$ sums up to $0$ unless $k$ is divisible by $M$. Note that in our alternative notation, the unit sector
corresponds to the ``root of unity'' $\zeta =1$, of primitive order  $m(\zeta)=1$ and of level $r(\zeta)=M$. Replacing also $k':=k/M$ by $k$,
we obtain 
\[ (1-q^r)\, e^{\textstyle \sum_{k\neq 0} \Psi^{rk}[(E^{(rk)}(q)-E^{(rk)}(1))/k(1-q)]} .\]
Notice that this answer coincides with \[
\Psi^r \left( (1-q)\, e^{\textstyle \sum_{k\neq 0} \Psi^k(E^{(rk)}(q)-E^{(rk)}(1))/k(1-q^k)} \right)
,\]
i.e. with $\Psi^r(\vv_r)$, where $\vv_r$ is the $r$th component of the dilaton vector described in Theorem 2.


\begin{thebibliography}{10000}
  
 

  


  

\bibitem{COGP} P. Candelas, X.C. de la Ossa, P.S. Green, L. Parkes. {\em A pair of Calabi–Yau manifolds as an exactly soluble superconformal theory.} Nucl. Phys. B359 (1991), 21.
  

  
\bibitem{Co} T. Coates. {\em Riemann--Roch theorems in Gromov--Witten theory.} PhD thesis, 2003, available at http://math.harvard.edu/~tomc/thesis.pdf 

\bibitem{CoGi} T. Coates, A. Givental. {\em Quantum Riemann--Roch, Lefschetz and Serre.} Ann. of Math. (2), 165 (2007), 15--53.

\bibitem{CGL} T. Coates, A. Givental. {\em Quantum cobordisms and formal group laws.} The unity of mathematics, 155--171, Progr. Math., 244, Birkh\"auser Boston, Boston, MA, 2006.




  
  
  

  
\bibitem{GiH} A. Givental. {\em Homological geometry I. Projective hypersurfaces.} Selecta Math. (New Series) 1 (1995), 325--345. 
 


  

\bibitem{GiQ} A. Givental. {\em Gromov-Witten invariants and quantization of quadratic Hamiltonians.} Mosc. Math. J. 1 (2001), no. 4, 551--568, 645 (English, with English and Russian summaries).

\bibitem{GiF} A. Givental. {\em Symplectic geometry of Frobenius structures.} Frobenius manifolds, Aspects Math., E36, Vieweg, Wiesbaden, 2004, pp. 91--112.


  
\bibitem{GiL} A. Givental, Y.-P. Lee. {\em Quantum K-theory on flag manifolds, finite difference Toda lattices and quantum groups.} Invent. Math. 151, 193--219, 2003.



  


\bibitem{KK} A. Kabanov, T. Kimura. {\em A change of coordinates on the large 
phase space of quantum cohomology.} Comm. Math. Phys. 217 (2001),
no. 1, 107--126, arXiv:math/9907096





\bibitem{KPSZ} P. Koroteev, P. Pushkar, A. Smirnov, A. Zeitlin. {\em Quantum K-theory of quiver varieties and many-body systems.} Preprint, 2017, 29 pp., arXiv:1705.10419 

  
  
\bibitem{YPLee} Y.-P. Lee. {\em Quantum K-theory I. Foundations.} Duke Math. J. 121 (2004), no. 3, 389--424.



  

\bibitem{ToK} V. Tonita. {\em A virtual Kawasaki Riemann–Roch formula.} Paciﬁc J. Math. 268 (2014), no. 1, 249--255. arXiv:1110.3916. 

\bibitem{ToT} V. Tonita. {\em Twisted orbifold Gromov--Witten invariants.} Nagoya Math. J. 213 (2014), 141--187, arXiv:1202.4778

\bibitem{ToN} V. Tonita. {\em Twisted K-theoretic Gromov--Witten invariants.} Preprint, 26 pp., arXiv: 1508.05976


  
\bibitem{Ts} H.-H. Tseng. {\em Orbifold quantum Riemann-Roch, Lefschetz and Serre.} Geom. Top. 14 (2010), 1--81.  
  

\end{thebibliography}
\enddocument